\documentclass[12pt,reqno]{amsart}

\usepackage[textheight=22cm,textwidth=16.6cm,headsep=1cm,centering]{geometry}
\usepackage{color}
\usepackage{esint,amssymb}
\usepackage{graphicx}
\usepackage{tikz}

\usepackage{hyperref}

\def\XXint#1#2#3{{\setbox0=\hbox{$#1{#2#3}{\int}$ }
\vcenter{\hbox{$#2#3$ }}\kern-.6\wd0}}

\newcommand{\R}{\mathbb{R}}

\newcommand{\ep}{\varepsilon}

\begin{document}

\title[Regularity of stable solutions to reaction-diffusion elliptic equations]
{Regularity of stable solutions to reaction-diffusion elliptic equations}

\author[Xavier Cabr\'e]{Xavier Cabr\'e}
\address{﻿X.C.\textsuperscript{1,2} ---
\textsuperscript{1}ICREA, Pg.\ Lluis Companys 23, 08010 Barcelona, Spain \& 
\textsuperscript{2}Universitat Polit\`ecnica de Catalunya, Departament de Matem\`{a}tiques and IMTech, 
Diagonal 647, 08028 Barcelona, Spain.
}
\email{xavier.cabre@upc.edu}

\begin{abstract}
The boundedness of stable solutions to semilinear (or reaction-diffusion) elliptic PDEs has been studied since the 1970's. In dimensions 10 and higher, there exist stable energy solutions which are unbounded (or singular). This note describes, for non-expert readers, a recent work in collaboration with Figalli, Ros-Oton, and Serra, where we prove that stable solutions are smooth up to the optimal dimension 9. This answers to an open problem posed by Brezis in the mid-nineties concerning the regularity of extremal solutions to Gelfand-type problems.

We also describe, briefly, a famous analogue question in differential geometry: the regularity of stable minimal surfaces.
\end{abstract}

\thanks{The author was supported by MINECO grant MTM2017-84214-C2-1-P}

\maketitle

\section{Hilbert's 19th problem and the principle of least action}

In many physical phenomena and geometric problems, observable states try to minimize a certain functional. When we describe the possible states by functions $u$ of one or several real variables, the functional is a real valued function $A$ acting on such functions. In classical mechanics, $A$ is called the {\it action} and is given by the integral of a Lagrangian. A simple example is the motion of a particle under gravitation, in which its position is given by $u=u(x)$ (where $x=t\in\R$ is time) and the action is the difference of kinetic and potential energies. In geometry, two important examples are geodesics (curves in a Riemannian manifold that are critical points of the length functional) and minimal surfaces (hypersurfaces of Euclidean space that are critical points of the area functional).

Hilbert's 19th problem asks whether minimizers of elliptic functionals are always analytic. When the functional is given by $A(u)=\int_{\Omega} L\left( \nabla u(x)\right) dx$ for some domain $\Omega\subset\R^n$ and convex function $L:\R^n\to\R$ (here $u:\Omega\subset\R^n\to\R$), the problem was solved independently in the late 1950s by Ennio De Giorgi and John Forbes Nash, Jr. Our work \cite{CFRS} takes up on the same question for the Lagrangian $L( \nabla u, u)=\frac{1}{2}|\nabla u|^2- F(u)$, which depends also on the variable~$u$.

The {\it principle of least action} in mechanics states that observable states should be not only critical points of the action but absolute minimizers of it. In the previous setting, when $L=L(\nabla u)$ is a convex function defined in $\R^n$, we have that $A$ is convex. 
As a consequence, any critical point of $A$ (if it exists) is an absolute minimizer. Thus, the principle of least action is, here, fulfilled. 
However, the principle is violated in some real life situations, described next, in which we observe states that are only local minimizers (minimizers among small perturbations), even in situations when an absolute minimizer exists. Still, the observed state being a local minimizer, it is therefore a {\it stable state}, in the sense that the second variation of the functional is nonnegative definite when computed at such state.

In these minimization problems, the competitors among which one looks for critical points are functions, or surfaces, with prescribed given boundary values ---the end points of the trajectory of a mechanical particle, or of a geodesic, or a given wire from which soap films or (minimal) surfaces are spanned. The key point is that the functional in many of these variational problems is not convex. Thus, it may admit critical points which are not absolute minimizers (but are only local minimizers) and even unstable critical points (which, consequently, are not even local minimizers).

\section{Stable minimal surfaces}

An instructive example is that of catenoids: a soap film or minimal surface formed between two coaxial parallel circular rings. In the nice paper ``{\it In situ observation of a soap-film catenoid---a simple educational physics experiment}'' by Ito and Sato \cite{IS}, catenoids are experimentally produced in a lab and videotaped while the distance between the two circular wires is continuosly increased. Note that, besides catenoids, there always exists another critical point of the area: the two flat disks spanned by the wires. For each small enough distance, two catenoids exist and the one with a thicker neck is an absolute minimizer of the area (clearly the disks have much larger area). As the wires separate, there is a distance $h_0$ at which both states (the thick-neck catenoid and the two disks) have the same area. Right after it, the two disks become the absolute minimizer, while the catenoid is only a local minimizer. Still, for an interval of distances $h>h_0$ the videotaped surface is the catenoid ---not the absolute minimizing disks. Such catenoid is a {\it stable minimal surface} ---stability understood as defined in the previous section.  Up to these distances, the unobserved thiner neck catenoid always existed and was an unstable minimal surface (this fits with the idea that a functional with two local minima should have a third unstable critical point). Finally, there is a second larger distance $h_c$ at which the local minimizer (the thick catenoid) and the unstable critical point (the thin one) get together to produce an inflection point. Right after it, the two-disks is the only critical point. In the experiment \cite{IS}, the unstable catenoid is photographed for a short instant, right before the distance $h_c$. Quickly after such instant, the thin neck collapses and the catenoid film succeeds to transform itself into the two disks.

The regularity theory of minimal surfaces has been the source of many important progresses in the area of PDEs.
In the 1960s the Italian school proved that the Simons cone
$$
x_1^2+\ldots+x_m^2 = x_{m+1}^2+\ldots+x_{2m}^2 
$$
is an absolute minimizer of area (for its own boundary values in any ball)
 if $2m \geq 8$, while it is not even stable in dimensions 2, 4, and~6.\footnote{This different behavior can be roughly understood noticing that the Jacobian for area in spherical coordinates, $r^{2m-2} dr$, becomes smaller at the origin as the dimension $2m$ increases. Note that for $2m=2$, the minimizer clearly avoids the origin:  for the boundary values of the cone, it is given by two parallel lines (and not by the ``cross'' passing through the origin).}
Thus, minimizing minimal surfaces of dimension $n$ may have (conical) singularities when $n\geq 7$. At the same time, a sequence of outstanding~contributions by different authors (J. Simons' being a prominent one) established that 
$n$-dimensional (absolute) minimizing minimal surfaces in $\R^{n+1}$ are always smooth when~$n\leq 6$. 
 
It is a long-standing open problem to extend this regularity result to the larger class of {\it stable minimal surfaces}. It is only known to be true for surfaces of dimension 2 or 3. See the survey \cite{CM} for more details on these issues.

\section{Stable solutions to reaction-diffusion elliptic equations}

The paper \cite{CFRS} takes on the analogue question (the regularity of {\it stable solutions}) for equations of the form $-\Delta u = f(u)$, where $\Delta$ is the Laplacian. They are called semilinear or reaction-diffusion elliptic equations and arise in many physical and biological situations. In the following combustion problem, 
a similar phenomenon to that of  catenoids occurs. It concerns the thermal self-ignition of a chemically active mixture of gases in a container. The model was introduced by Frank-Kamenetskii in the 1930s but became popular within the mathematical community when Barenblatt wrote chapter~15 of the volume \cite{Gelf}, edited by Gelfand in 1963. Here $x\in \Omega\subset \R^n$ denotes points in the container~$\Omega$ and $u=u(x)$ is the temperature at the point. The action functional is the difference of kinetic and potential energies,
$$
A(u)=\int_{\Omega}\Bigl(\frac{1}{2}|\nabla u(x)|^2-F(u(x))\Bigr)\,dx,
$$
where $F:\R\to\R$ is a given function, which Barenblatt chose to be $F(u)=\lambda e^u$, with $\lambda$ a positive constant, from Arrhenius law in chemical kinetics. For convenience, we will impose vanishing boundary conditions: the temperature is kept at $u=0$ on the boundary $\partial\Omega$. Making a first variation $u+\ep v$ and integrating by parts, one easily sees that critical points of $A$ satisfy the {\it reaction-diffusion equation}
\begin{equation}
\label{eq:PDE}
-\Delta u=f(u) \quad \text{ in }\Omega\subset \R^n,
\end{equation}
where $f=F'$. In the case (among others) of the so called {\it Gelfand problem}, $-\Delta u=\lambda e^u$, the situation is similar to the one of catenoids. For a certain range $\lambda\in (0,\lambda^*)$ of parameters, there exists a stable solution $u_\lambda$ ---that is, a solution at which the functional $A$ has a nonnegative definite second variation. Such stable solution is not an absolute minimizer since $A$ is unbounded by below (note that the potential density $F(u)=e^u$ grows faster at infinity than the quadratic kinetic density $|\nabla u(x)|^2$). For some parameters $\lambda \in (0,\lambda^*)$, there might also exist (this will depend on the container~$\Omega$) unstable solutions of the same problem. For $\lambda=\lambda^*$, the limit of the functions $u_\lambda$ is an $L^1$ weak stable solution, called the extremal solution. When $\lambda> \lambda^*$, no solution exists ---in the same way that catenoids did not exist for distances $h$ between the wires larger than $h_c$.

To better understand the problem, let us consider the nonlinear heat equation 
\begin{equation}
\label{eq:heat}
v_t-\Delta v=f(v), 
\end{equation}
where $v=v(x,t)$ and $t$ is time. Now, a stable solution $u=u(x)$ of \eqref{eq:PDE} can be understood as a stationary solution of \eqref{eq:heat} which is stable in the sense of Lyapunov ---note that a simple computation shows that the action functional $A(v(\cdot,t))$ is non-increasing in the time $t$. The problem is nonlinear due to the sources of heat, $f(v(x,t))$ or $f(u(x))$: the production of heat depends nonlinearly on the actual temperature. As described in \cite{GM,GMN}, equation~\eqref{eq:heat}
describes the evolution of an initially uniform temperature $v(\cdot,0)\equiv 0$ which diffuses in space and increases in the container due to the heat release given by the reaction term $f(v)$ ---note that in Gelfand's problem the initial heat source $\lambda f(0)= \lambda e^0=\lambda$ is already positive. The parameters $\lambda$ for which there exists a stable solution  of \eqref{eq:PDE} correspond to ignition failure (the reactive component undergoes partial oxidation and results in establishing a stationary temperature profile equal to the stable solution). Instead, $\lambda>\lambda^*$ (when there exists no stationary solution) means successful auto-ignition in the combustion process.

Since we will turn now to regularity issues, let us recall that Fourier invented his omnipresent Fourier series to understand the linear heat equation. On the other hand, the regularity theory for the stationary linear Poisson equation $-\Delta u=g(x)$ is at the center of PDE theory and also propitiated the development of many tools, such as the theory of singular integrals in harmonic analysis. For instance, it is well known the Lebesgue-integrability requirements for $g=g(x)$ needed to guarantee the boundedness of the potential function~$u$. This is relevant since, for our nonlinear equation~\eqref{eq:PDE}, $-\Delta u= f(u)$, the obstruction for regularity is the possibility that $u$ becomes unbounded somewhere, that is, $u$ blows-up at some points (and hence $u$ is singular). As we will see next, this singular behavior can be produced by the strength of some reaction terms $f(u)$. A technical detail for experts is that, since our problem is variational, in what follows we consider only (singular) solutions for which each term of the action functional is integrable (that is, energy solutions).

When $n\geq 3$, $\Omega=B_1$ is the unit ball, 
\begin{equation*}\label{eq:n10}
u=\log\frac{1}{|x|^2}, \quad \text{and}\quad f(u)=2(n-2)e^u,
\end{equation*}
a simple computation shows that we are in the presence of a singular solution of \eqref{eq:PDE} vanishing on $\partial B_1$. As the Simons cone in minimal surfaces theory, this explicit solution turns out to be stable in high dimensions, precisely when $n\geq 10$. On the other hand, in the 1970s Crandall and Rabinowitz~\cite{CR} established that if
\begin{equation*}
f(u)=e^u \quad\text{ or }\quad f(u)=(1+u)^p \, \text{ with } p>1,
\end{equation*}
then stable solutions in any smooth bounded domain $\Omega$ are bounded (and hence smooth and analytic, by classical elliptic regularity theory) when $n \leq 9$.
This result was the main reason for Haim Brezis to raise the following question in the 1990s (which we cite almost literally from a later reference):
 
\vspace{2mm}

\noindent
Brezis \cite[Open problem 1]{Br}: {\it Is there something ``sacred'' about dimension 10? More precisely, is it possible in ``low'' dimensions to construct some $f$ (and some $\Omega$) for which a singular stable solution exists? Alternatively, can one prove in ``low'' dimensions that every stable solution is smooth for every $f$ and every $\Omega$? }

\vspace{2mm}

\noindent
Other open questions on stable solutions were posed by Brezis and V\'azquez \cite{BV}. 

The last twenty five years have produced a large literature on Gelfand-type problems. See the monograph~\cite{Dup} for an extensive list of results and references. For a certain type of nonlinearities~$f$, some of these works are related to micro-electro-mechanical systems (MEMS); see~\cite{EGG}. 

The main developments proving that stable solutions to \eqref{eq:PDE} are smooth (no matter what the nonlinearity $f$ is) were made

\begin{itemize}
\item by Nedev \cite{Ned00} in 2000, when $n \leq 3$ (and $f$ is convex);
\item by Cabr\'e and Capella  \cite{CC06} in 2006, when $\Omega=B_1$ ($u$ is radially symmetric) and $n\leq 9$;
\item by Cabr\'e \cite{C10} in 2010, when $n\leq 4$ (and $\Omega$ is convex).
\end{itemize}
Note that the 2006 result in the radially symmetric case, \cite{CC06}, accomplished the optimal dimension $n\leq 9$ for every nonlinearity $f$. This gave hope for the result to be true also in the general nonradial case, though no certainty was assured---note that Brezis' statement above leaves both the affirmative and negative answers  as possible ones. Since 2010, after~\cite{C10}, the regularity result was only known up to dimension $n=4$. Two attempts in higher dimensions (recorded in~\cite{CFRS}) gave only very partial answers.

The work \cite{CFRS} finally solves the open problem, by establishing the regularity of stable solutions to \eqref{eq:PDE} in the interior of any open set $\Omega$ in the optimal dimensions $n\leq 9$ under the only requirement for the nonlinearity $f$ to be nonnegative. Furthermore, adding the vanishing boundary condition $u=0$ on $\partial\Omega$, the article proves regularity up to the boundary when $\Omega$ is of class $C^3$ and $n\leq 9$, assuming now $f$ to be nonnegative, nondecreasing, and convex. Both results come along with new universal H\"older-continuity estimates which have a very weak norm (the $L^1$-norm) of the solution on their right-hand sides. They read, respectively, as
$$
\|u\| _{C^\alpha(\overline{B}_{1/2})}\leq C\|u\|_{L^1(B_1)}
\quad\text{ and }\quad
\|u\| _{C^\alpha(\overline\Omega)} \leq C_\Omega \, \|u\|_{L^1(\Omega)},
$$
where $\alpha>0$ and $C$ are dimensional constants, while $C_\Omega$ depends only on $\Omega$. These estimates are rather surprising (because of their universality) for a nonlinear problem, specially since they make no reference to the reaction nonlinearity $f$. The stability of the solution $u$ is crucial for their validity.  For the expert reader, \cite{CFRS} also establishes another open problem from \cite{BV}: an apriori $H^1=W^{1,2}$ estimate for stable solutions in all dimensions $n$. 

If the nonnegativeness of $f$ is a needed requirement for interior regularity remains as an open question. It is only known to be unnecessary for  $n\leq 4$, as well as for $n\leq 9$ in the radial case.

The proofs in the article are too technical to be described here. Let us only say that a key point is to use the stability property under two different types of small perturbations of the solution $u$ (one in the radial direction, the other in the normal direction to the level sets):
$$
u\left( x+ \ep |x|^{(2-n)/2}\zeta (x) x \right)
\quad\text{ and }\quad
u\Big( x+ \ep \eta (x)\frac{\nabla u (x)}{|\nabla u(x)|}  \Big),
$$
where $\zeta$ and $\eta$ are cut-off functions.

After \cite{CFRS}, an analogue result for equations involving the $p$-Laplacian has been proved by Cabr\'e, Miraglio, and Sanch\'on~\cite{CMS}. It is optimal in terms of dimensions for $p>2$, but not for $p<2$ ---this case remains as an open problem. On the other hand, for the recently very active area of fractional Laplacians, an optimal result for $(-\Delta)^s u=f(u)$ is largely open ---even in the radial case. The optimal dimensions for regularity have only been accomplished in a 2014 work of Ros-Oton~\cite{R-O} for the Gelfand nonlinearity $f(u)=\lambda e^u$ in symmetric convex domains ---but for any fraction $s\in (0,1)$ of the Laplacian.

\end{document}